\input amstex
\input amsppt.sty   
\hsize 30pc
\vsize 47pc
\def\nmb#1#2{#2}         
\def\cit#1#2{\ifx#1!\cite{#2}\else#2\fi} 
\def\totoc{}             
\def\idx{}               
\def\ign#1{}             

\redefine\o{\circ}

\define\al{\alpha}

\define\ga{\gamma}
\define\de{\delta}
\define\ep{\varepsilon}

\define\la{\lambda}

\define\si{\sigma}

\define\ps{\psi}

\define\De{\Delta}

\define\La{\Lambda}

\redefine\i{^{-1}}
\define\row#1#2#3{#1_{#2},\ldots,#1_{#3}}
\define\x{\times}
\def\today{\ifcase\month\or
 January\or February\or March\or April\or May\or June\or
 July\or August\or September\or October\or November\or December\fi
 \space\number\day, \number\year}
\topmatter
\title Choosing roots of polynomials smoothly
\endtitle
\author Dmitri Alekseevsky\\
Andreas Kriegl\\ Mark Losik\\ Peter W\. Michor  \endauthor
\leftheadtext{D\.V\. Alekseevsky, A\. Kriegl, M\. Losik, P\.W\. Michor}
\affil
Erwin Schr\"odinger International Institute of Mathematical Physics, 
Wien, Austria
\endaffil
\address 	 D\. V\. Alekseevsky: Center `Sophus Lie', 
Krasnokazarmennaya 6, 111250 Moscow, Russia
\endaddress
\address
A\. Kriegl, P\. W\. Michor: Institut f\"ur Mathematik, Universit\"at Wien,
Strudlhofgasse 4, A-1090 Wien, Austria
\endaddress
\email Andreas.Kriegl\@univie.ac.at, Peter.Michor\@esi.ac.at \endemail
\address M\. Losik: Saratov State University, ul\. Astrakhanskaya, 83, 
410026 
Saratov, Russia \endaddress
\email losik\@scnit.saratov.su\endemail
\date {\today} \enddate
\thanks 
Supported by `Fonds zur F\"orderung der wissenschaftlichen 
Forschung, Projekt P~10037~PHY'.
\endthanks
\abstract We clarify the question whether for a smooth curve of 
polynomials one can choose the roots smoothly and related questions.
Applications to perturbation theory of operators are given.
\endabstract
\endtopmatter

\document

\head Table of contents \endhead
\noindent 1. Introduction\leaders \hbox to 
1em{\hss .\hss }\hfill {\eightrm 1}\par  
\noindent 2. Choosing differentiable square and cubic roots 
\leaders \hbox to 1em{\hss .\hss }\hfill {\eightrm 2}\par  
\noindent 3. Choosing local roots of real polynomials smoothly 
\leaders \hbox to 1em{\hss .\hss }\hfill {\eightrm 7}\par  
\noindent 4. Choosing global roots of polynomials differentiably 
\leaders \hbox to 1em{\hss .\hss }\hfill {\eightrm 11}\par  
\noindent 5. The real analytic case \leaders \hbox to 
1em{\hss .\hss }\hfill {\eightrm 15}\par  
\noindent 6. Choosing roots of complex polynomials \leaders \hbox to 
1em{\hss .\hss }\hfill {\eightrm 16}\par  
\noindent 7. Choosing eigenvalues and eigenvectors of matrices 
smoothly \leaders \hbox to 1em{\hss .\hss }\hfill {\eightrm 18}\par  

\head \totoc \nmb0{1}. Introduction\endhead 

We consider the following problem.
Let
$$
P(t)=x^n-a_1(t)x^{n-1}+\dots+(-1)^na_n(t) \tag1
$$
be a polynomial with all roots real, smoothly parameterized by $t$ near 
0 in $\Bbb R$. 
Can we find $n$ smooth functions $x_1(t),\dots,x_n(t)$ of the 
parameter $t$ defined near 0, which are the roots of $P(t)$ for each $t$? 
We can reduce the problem to $a_1=0$, replacing the variable $x$ by 
with the variable $y=x-a_1(t)/n$.
We will say that the curve (1) is smoothly solvable near $t=0$ if 
such smooth roots $x_i(t)$ exist.

We describe an algorithm which in the smooth and in the holomorphic 
case sometimes allows to solve this problem.
The main results are:
If all roots are real then they can always be chosen differentiable, 
but in general not $C^1$ for degree $n\ge3$; and in degree 2 they can 
be chosen twice differentiable but in general not $C^2$.
If they are arranged in increasing order, they depend continuously on 
the coefficients of the polynomial, and if moreover no two of them 
meet of infinite order in the parameter, then they can be chosen 
smoothly.
We also apply these results to obtain a smooth 1-parameter 
perturbation theorem for selfadjoint operators with compact resolvent 
under the condition, that no pair of eigenvalues meets of infinite 
order.

We thank C\. Fefferman who found a mistake in a first version of 
\nmb!{2.4}, to M\. and Th\. Hoffmann-Ostenhof for their interest 
and hints, and to Jerry Kazdan for arranging \cit!{15}.

\head\totoc\nmb0{2}. Choosing differentiable square and cubic roots \endhead

\proclaim{\nmb.{2.1}. Proposition. The case $n=2$}
Let $P(t)(x)=x^2-f(t)$ for a function $f$ defined on an open 
interval, such that $f(t)\geq 0$ for all $t$.

If $f$ is smooth and is nowhere flat of infinite order, 
then smooth solutions $x$ exist.

If $f$ is $C^2$ then $C^1$-solutions exist.

If $f$ is $C^4$ then twice differentiable solutions exist.
\endproclaim

\demo{Proof} Suppose that $f$ is smooth.
If $f(t_0)>0$ then we have obvious local smooth solutions 
$\pm\sqrt{f(t)}$.
If $f(t_0)=0$ we have to find a smooth function $x$ such that $f=x^2$, 
a smooth square root of $f$.
If $f$ is not flat at $t_0$ then the 
first nonzero derivative at $t_0$ has even order $2m$ and is positive, 
and $f(t)=(t-t_0)^{2m}f_{2m}(t)$, where 
$f_{2m}(t):=\int_0^1\frac{(1-r)^{2m-1}}{(2m-1)!}f^{(2m)}(t_0+r(t-t_0))\,dr$ 
gives a smooth function and 
$f_{2m}(t_0)=\frac1{(2m)!}f^{(2m)}(t_0)>0$.
Then $x(t):= (t-t_0)^m\sqrt {f_{2m}(t)}$ is a local smooth solution. 
One can piece together these local solutions, changing sign at all 
points where the first non-vanishing derivative of $f$ is of order 
$2m$ with $m$ odd.
These points are discrete. 

Let us consider now a function $f\ge 0$ of class $C^2$. 
We claim that then $x^2=f(t)$ admits a $C^1$-root $x(t)$, 
globally in $t$. 
We consider a fixed $t_0$.
If $f(t_0)>0$ then there is locally even a $C^2$-solution 
$x_\pm(t)=\pm\sqrt{f(t)}$.  
If $f(t_0)=0$ then $f(t)=(t-t_0)^2h(t)$ 
where $h\ge0$ is continuous and $C^2$ off $t_0$ with 
$h(t_0)=\tfrac12f''(t_0)$.
If $h(t_0)>0$ then $x_{\pm}(t)=\pm (t-t_0)\sqrt{h(t)}$ is 
$C^2$ off $t_0$, and
$$
x_\pm'(t_0)=\lim_{t\to t_0} \frac{x_{\pm}(t)-x_{\pm}(t_0)}{t-t_0} =
     \lim_{t\to t_0} \pm\sqrt{h(t)} =  \pm\sqrt{h(t_0)} = 
     \pm\sqrt{\tfrac12f''(t_0)}.
$$
If $h(t_0)=0$ then we choose $x_{\pm}(t_0)=0$, and any choice of the 
roots is then differentiable at $t_0$ with derivative $0$, 
by the same calculation.  

One can piece together 
these local solutions: At zeros $t$ of $f$ where  $f''(t)>0$ we have 
to pass through 0, but where $f''(t)=0$ the choice of the root does 
not matter.
The set $\{t: f(t)=f''(t)=0\}$ is closed, so its 
complement is a union of open intervals.
Choose a point in each of 
these intervals where $f(t)>0$ and start there with the positive 
root, changing signs at points where $f(t)=0\ne f''(t)$: these points 
do not accumulate in the intervals.
Then we get 
a differentiable choice of a root $x(t)$ on each of this open intervals 
which extends to a global differentiable root which is 0 on $\{t: 
f(t)=f''(t)=0\}$.  

Note that we have 
$$
x'(t)=\cases \frac{f'(t)}{2x(t)} & \text{ if }f(t)> 0\\
                   \pm\sqrt{f''(t)/2} & \text{ if }f(t)=0\endcases
$$
In points $t_0$ with $f(t_0)>0$ the solution $x$ is 
$C^2$; locally around points $t_0$ with 
$f(t_0)=0$ and $f''(t_0)>0$ the root $x$ is $C^1$ since 
for $t\ne t_0$ near $t_0$ we have $f(t)>0$ and $f'(t)\ne 0$, so by 
l'Hospital we get
$$\align
\lim_{t\to t_0}x'(t)^2 &= \lim_{t\to t_0}\frac{f'(t)^2}{4f(t)} 
=  \lim_{t\to t_0}\frac{2f'(t)f''(t)}{4f'(t)} = \frac{f''(t_0)}2 = 
x'(t_0)^2,  
\endalign$$
and since the choice of signs was coherent, $x'$ is continuous at $t_0$;
if $f''(t_0)=0$ then $x'(t_0)=0$ and $x'(t)\to 0$ 
for $t\to t_0$ for both expressions, by lemma \nmb!{2.2} below. 
Thus $x$ is $C^1$.

If moreover $f\ge 0$ is $C^4$, then the solution $x$ from above may 
be modified to be twice differentiable.
Near points $t_0$ with $f(t_0)>0$ any continuous root 
$t\mapsto x_\pm(t)=\pm\sqrt{f(t)}$ is even $C^4$. 
Near points $t_0$ with $f(t_0)=f'(t_0)=0$ we 
have $f(t)=(t-t_0)^2h(t)$ where $h\ge0$ is $C^2$.
We may choose 
a $C^1$-root $z$ with $z^2=h$ by the arguments above, and then 
$x(t):=(t-t_0)z(t)$ is twice differentiable at 
$t_0$ since we have 
$$\align
\frac{x'(t)-x'(t_0)}{t-t_0} &= \frac{z(t)+(t-t_0)z'(t)-z(t_0)}{t-t_0} \\
& =z'(t) + \frac{z(t)-z(t_0)}{t-t_0} \to 2z'(t_0) = 
     \pm\sqrt{\frac{f^{(4)}(t_0)}{4!}}.
\endalign$$ 
If $f(t_0)=f''(t_0)=f^{(4)}(t_0)=0$ then any $C^1$ choice of the roots is 
twice differentiable at $t_0$, in particular $x(t)=|t-t_0|z(t)$.

Now we can piece together this solutions similarly as above.
Let $y$ be a global $C^1$-root of $f$, chosen as above changing
sign only at points $t$ with $f(t)=0<f''(0)$.
We put $x(t)=\ep(t) y(t)$, where $\ep(t)\in\{\pm1\}$ will be chosen 
later. 
The set $\{t:f(t)=f''(t)=f^{(4)}(t)=0\}$ has a countable union of 
open intervals as complements.
In each of these intervals choose 
a point $t_0$ with $f(t_0)>0$, near which $y$ is $C^4$.
Now let
$\ep(t_0)=1$, and let $\ep$ change sign 
exactly at points with $f(t)=f''(t)=0$ but $f^{(4)}(t)>0$.
These 
points do not accumulate inside the interval.
Then $x$ is twice 
differentiable. 
\qed\enddemo  
                                     
\proclaim{\nmb.{2.2}. Lemma}
Let $f\ge 0$ be a $C^2$-function with $f(t_0)=0$, then 
for all $t\in \Bbb R$ we have 
$$
f'(t)^2\le 2 f(t)\;\max\{f''(t_0+r(t-t_0)): 0\le r\le 2\}. \tag1
$$
\endproclaim

\demo{Proof}
If $f(t)=0$ then $f'(t)=0$ so \thetag 1 holds. 
We use the Taylor formula
$$
f(t+s) = f(t) + f'(t)s + \int_0^1(1-r)f''(t+rs)\,dr\;s^2 \tag2
$$
In particular we get (replacing $t$ by $t_0$ and then $t_0+s$ by $t$)
$$\align
f(t) &= 0 + 0 +  \int_0^1(1-r)f''(t_0+r(t-t_0))\,dr\;(t-t_0)^2 \tag3 \\
&\le \frac{(t-t_0)^2}2 \max\{f''(t_0+r(t-t_0)): 0\le r\le 2\} 
\endalign$$
Now in \thetag2 we replace $s$ by $-\ep s$ (where 
$\ep=\operatorname{sign}(f'(t))$) to obtain 
$$
0\le f(t-\ep s) = f(t) - |f'(t)|s 
     +  \int_0^1(1-r)f''(t-\ep rs)\,dr\;s^2 \tag4
$$
Let us assume $t\ge t_0$ and then choose (using \thetag3)
$$
s(t):= \sqrt{\frac{2f(t)}{\max\{f''(t_0+r(t-t_0)): 
0\le r\le 2\}}}\le t-t_0.
$$
Note that we may assume $f(t)>0$, then $s(t)$ is well defined and 
$s(t)>0$.
This choice of $s$ in \thetag 4 gives
$$\align
|f'(t)| &\le \frac1{s(t)}\left(f(t)+\frac{s(t)^2}2\max\{f''(t-\ep rs(t)): 
     0\le r\le 1\}\right) \\
&\le \frac1{s(t)}\left(f(t)+\frac{s(t)^2}2\max\{f''(t-r(t-t_0)): 
     -1\le r\le 1\}\right) \\
&= \frac{2f(t)}{s(t)} 
     = \sqrt{2 f(t)\;\max\{f''(t_0+r(t-t_0)): 0\le r\le 2\}}
\endalign$$
which proves \thetag1 for $t\ge t_0$.
Since the assertion is 
symmetric it then holds for all $t$.
\qed\enddemo

\subhead\nmb.{2.3}. Examples \endsubhead
If $f\ge0$ is only $C^1$, then there may not
exist a differentiable root of $x^2=f(t)$, as the following example 
shows: 
$x^2=f(t):= t^2\sin^2(\log t)$ is $C^1$, but $\pm t\sin(\log t)$ is not 
differentiable at 0.

If $f\ge0$ is twice differentiable there may not exist a $C^1$-root:
$x^2=f(t)=t^4\sin^2(\tfrac1t)$ is twice differentiable, but 
$\pm t^2\sin(\tfrac1t)$ is differentiable, but not $C^1$. 

If $f\ge0$ is only $C^3$, then there may not
exist a twice differentiable root of $x^2=f(t)$, as the following example 
shows: 
$x^2=f(t):= t^4\sin^2(\log t)$ is $C^3$, but $\pm t^2\sin(\log t)$ is 
only $C^1$ and not twice differentiable.  

\subhead {\nmb.{2.4}. Example} \endsubhead 
If $f(t)\ge0$ is smooth but flat at $0$, in general our problem has no 
$C^2$-root as 
the following example shows, which is an application of the general 
curve lemma 4.2.15 in \cit!{5}: Let $h:\Bbb R\to [0,1]$ be smooth 
with $h(t)=1$ for $t\ge0$ and $h(t)=0$ for $t\le-1$.
Then the 
function
$$\align
f(t) :&= \sum_{n=1}^\infty h_n(t-t_n). 
     \left(\frac{2n}{2^n}(t-t_n)^2+\frac1{4^n}\right),\qquad \text{ where}\\
h_n(t) :&= h\left(n^2\left(\frac1{n.2^{n+1}}+t\right)\right). 
     h\left(n^2\left(\frac1{n.2^{n+1}}-t\right)\right)\quad\text{ and}\\
t_n :&= \sum_{k=1}^{n-1} \left(\frac 2{k^2}+\frac2{k.2^{k+1}}\right)
     + \frac 1{n^2}+\frac1{n.2^{n+1}},
\endalign$$
is $\ge 0$ and is smooth: the sum consists of at most one summand for 
each $t$, 
and the derivatives of the summands converge uniformly to 0:
Note that $h_n(t)=1$ for $|t|\leq \frac1{n.2^{n+1}}$ and
$h_n(t)=0$ for $|t|\geq \frac1{n.2^{n+1}}+\frac1{n^2}$ hence
$h_n(t-t_n)\ne 0$ only for 
$r_n<t<r_{n+1}$, where 
$r_n:=\sum_{k=1}^{n-1} \left(\frac 2{k^2}+\frac2{k.2^{k+1}}\right)$.
Let $c_n(s):=\frac{2n}{2^n}s^2+\frac1{4^n}\geq 0$ and
$H_i := \sup\{|h^{(i)}(t)|:t\in\Bbb R\}$.
Then
$$
\align
n^2\sup\{|(h_n\cdot c_n)^{(k)}&(t)|:t\in\Bbb R\}
     = n^2\sup\{|(h_n\cdot c_n)^{(k)}(t)|:|t|\leq 
     \frac1{n.2^{n+1}}+\frac1{n^2}\} \\
&\leq n^2\sum_{i=0}^k \binom{k}{i} n^{2i} H_i
	\sup\{|c_n^{(k-i)}(t)|:|t|\leq \frac1{n.2^{n+1}}+\frac1{n^2}\}\\
&\leq \biggl(\sum_{i=0}^k \binom{k}{i} n^{2i+2} H_i \biggr)
	\sup\{|c_n^{(j)}(t)|:|t|\leq 2,\;j\leq k\}
\endalign
$$
and since $c_n$ is rapidly decreasing in $C^\infty(\Bbb R,\Bbb R)$ 
(i.e\. $\{p(n)\,c_n:n\in\Bbb N\}$ 
is bounded in $C^\infty(\Bbb R,\Bbb R)$ for each polynomial $p$) 
the right side of the inequality
above is bounded with respect to $n\in\Bbb N$ and hence the series 
$\sum_n h_n(\quad-t_n) c_n(\quad-t_n)$ converges
uniformly in each derivative, and thus represents an element of 
$f\in C^\infty(\Bbb R,\Bbb R)$.
Moreover we have
$$
f(t_n)=\frac1{4^n},\quad f'(t_n)=0,\quad f''(t_n)=\frac{2n}{2^{n-1}}.
$$ 
Let us assume that $f(t)=g(t)^2$ for $t$ 
near $\sup_nt_n<\infty$, where $g$ is twice differentiable.
Then 
$$\align
f'&=2gg' \\
f''&= 2gg''+2(g')^2\\
2ff''&= 4g^3 g'' + (f')^2 \\
2f(t_n)f''(t_n)&=4g(t_n)^3g''(t_n)+f'(t_n)^2
\endalign$$
thus 
$g''(t_n)= \pm 2n$, so $g$ cannot be 
$C^2$, and $g'$ cannot satisfy a local Lipschitz condition near 
$\lim t_n$.
Another similar example can be found in \nmb!{7.4} below.

According to \cit!{15}, some results of this section are 
contained in Frank Warners dissertation (around 1963, unpublished):
Non-negative smooth functions have $C^1$ square roots whose second 
derivatives exist everywhere.
If all zeros are of finite order there 
are smooth square roots.
However, there are examples not 
possessing a  $C^2$ square root.
Here is one: 
$$
f(t) = \sin^2(1/t)e^{-1/t} + e^{-2/t} \text{ for }t>0,\quad f(t) = 
0\text{ for }t\le0.
$$
This is a sum of two non-negative $C^\infty$ functions each of which 
has a $C^\infty$ square root.
 But the second derivative of the square 
root of $f$ is not continuous at the origin. 

In \cit!{6} Glaeser proved that a non-negative 
$C^2$-function on an open subset of 
$\Bbb R^n$ which vanishes of second order has a $C^1$ positive square 
root.
A smooth function $f\ge 0$ is constructed which is flat at 0 
such that the positive square root is not $C^2$. 
In \cit!{4} Dieudonn\'e gave shorter proofs of Glaeser's results.

\subhead\nmb.{2.5}. Example. The case $n\ge 3$ \endsubhead
We will construct a polynomial 
$$
P(t) = x^3 + a_2(t)x -a_3(t)
$$
with smooth coefficients $a_2$, $a_3$ with all roots real which does 
not admit $C^1$-roots.
Multiplying with other polynomials one then 
gets polynomials of all orders $n\ge 3$ which do not admit 
$C^1$-roots.

Suppose that $P$ admits $C^1$-roots $x_1$, $x_2$, $x_3$.
Then we have
$$\align
0 &= x_1 + x_2 + x_3 \\
a_2 &= x_1x_2 + x_2x_3 + x_3x_1 \\
a_3 &= x_1x_2x_3\\ 
0 &= \dot x_1 + \dot x_2 + \dot x_3 \\
\dot a_2 &= \dot x_1x_2+x_1\dot x_2 + \dot x_2x_3+x_2\dot x_3 + \dot x_3x_1+x_3\dot x_1 \\
\dot a_3 &= \dot x_1x_2x_3+x_1\dot x_2x_3+x_1x_2\dot x_3 
\endalign$$
We solve the linear system formed by the last three equations and get
$$\align
\dot x_1 &= \frac{\dot a_3-\dot a_2x_1}{(x_2-x_1)(x_3-x_1)} \\
\dot x_2 &= \frac{\dot a_3-\dot a_2x_2}{(x_3-x_2)(x_1-x_2)} \\
\dot x_3 &= \frac{\dot a_3-\dot a_2x_3}{(x_1-x_3)(x_2-x_3)}. \\
\endalign$$
We consider the continuous function
$$
b_3 := \dot x_1\dot x_2\dot x_3 
     = \frac{\dot a_3^3 + \dot a_2^2\dot a_3a_2 - \dot a_2^3a_3}
          {4a_2^3+27a_3^2}.
$$
For smooth functions $f$ and $\ep $ with $\ep^2\le 1$ we let
$$\align
u :&= -12 a_2 := f^2 \\
v :&= 108 a_3 :=\ep f^3.
\endalign$$
then all roots of $P$ are real since $a_2\le 0$ and 
$432(4a_2^3+27a_3^2) = v^2-u^3 = f^6(\ep^2-1)\le 0$.
We get then 
$$\align
11664 b_3 &=\frac{4\dot v^3-27 u\dot u^2\dot v+27 \dot u^3 v}{v^2-u^3}\\
&= 4\frac{f^3\dot\ep^3+9f^2\dot f\dot\ep^2\ep+27f\dot f^2\dot\ep(\ep^2-1) 
     +27 \dot f^3\ep (\ep^2-3)}{\ep^2-1}.
\endalign$$
Now we choose 
$$\align
f(t) :&= \sum_{n=1}^\infty h_n(t-t_n). 
     \left(\frac{n}{2^n}(t-t_n)\right),\\
\ep(t) :&= 1-\sum_{n=1}^\infty h_n(t-t_n).\frac1{8^n}, 
\endalign$$
where $h_n$ and $t_n$ are as in the beginning of \nmb!{2.4}.
Then $f(t_n)=0$, $\dot f(t_n)=\frac n{2^n}$, and 
$\ep(t_n)=\frac1{8^n}$, hence 
$$108 b_3(t_n) 
     = \frac{\dot f(t_n)^3\ep(t_n)(\ep(t_n)^2-3)}{\ep(t_n)^2-1} 
     \sim n^3 \to \infty$$
is unbounded on the convergent sequence $t_n$.
So the roots cannot be 
chosen locally Lipschitz, thus not $C^1$.

\head\totoc\nmb0{3}. Choosing local roots of real polynomials smoothly 
\endhead
 
\subhead {\nmb.{3.1}. Preliminaries}\endsubhead 
We recall some known facts on polynomials with real 
coefficients.
Let 
$$ 
P(x)=x^n-a_1x^{n-1}+\dots+(-1)^na_n
$$ 
be a polynomial with real coefficients $a_1,\dots,a_n$ and roots 
$x_1,\dots,x_n\in \Bbb C$. 
It is known that $a_i=\sigma_i(x_1,\dots,x_n)$, where 
$\sigma_i$ $(i=1,\dots,n)$ are the elementary symmetric 
functions in $n$ variables: 
$$
\si_i(x_1,\dots,x_n)=\sum_{1\le j_1<\dots<j_i\le n}x_{j_1}\dots x_{j_i}
$$ 
Denote by $s_i$ the Newton polynomials 
$\sum_{j=1}^nx_j^i$ which are related to the elementary symmetric 
function by 
$$ 
s_k-s_{k-1}\sigma_1+s_{k-2}\sigma_2+\dots+(-1)^{k-1}s_1\sigma_{k-1}+ 
(-1)^kk\sigma_k=0  \quad (k\leq n) \tag 1
$$ 
The corresponding mappings are related by a polynomial diffeomorphism 
$\ps^n$, given by \thetag1:
$$\align
\sigma^n:&=(\sigma_1,\dots,\sigma_n):\Bbb R^n\to \Bbb R^n\\
s^n:&=(s_1,\dots,s_n):\Bbb R^n\to \Bbb R^n \\
s^n:&=\ps^n\o\si^n
\endalign$$ 
Note that the Jacobian (the determinant of the derivative) 
of $s^n$ is $n!$ times the Vandermonde determinant:
$\det(ds^n(x))=n!\,\prod_{i>j}(x_i-x_j)
     =:n!\,\operatorname{Van}(x)$, and even the derivative itself 
$d(s^n)(x)$ equals the 
Vandermonde matrix up to factors $i$ in the $i$-th row.
We also have 
$\det(d(\ps^n)(x))=(-1)^{n(n+3)/2}n!=(-1)^{n(n-1)/2}n!$, and 
consequently 
$\det(d\si^n(x))=\prod_{i>j}(x_j-x_i)$.
We consider the so-called Bezoutiant  
$$ B:= \pmatrix 
s_0&s_1&\hdots&s_{n-1}\\ 
s_1&s_2&\hdots&s_n\\ 
\vdots&\vdots&&\vdots \\ 
s_{n-1}&s_n&\hdots&s_{2n-2} 
\endpmatrix. 
$$ 
Let $B_k$ be the minor formed by the first $k$ rows and columns of $B$. 
 From
$$ 
B_k(x)=\pmatrix 
1&1&\hdots&1\\ 
x_1&x_2&\hdots&x_n\\ 
\vdots&\vdots&&\vdots\\ 
x_1^{k-1}&x_2^{k-1}&\hdots&x_n^{k-1} 
\endpmatrix\cdot 
\pmatrix 
1&x_1&\hdots&x_1^{k-1}\\ 
1&x_2&\hdots&x_2^{k-1}\\ 
\vdots&\vdots&&\vdots\\ 
1&x_n&\hdots&x_n^{k-1} 
\endpmatrix $$ 
it follows that
$$ 
\Delta_k(x):=\det(B_k(x))
     =\sum_{i_1<i_2<\dots<i_k}(x_{i_1}-x_{i_2})^2\dots(x_{i_1}-x_{i_n})^2 
     \dots(x_{i_{k-1}}-x_{i_k})^2,\tag 2
$$ 
since for $n\x k$-matrices $A$ one has 
$\det(AA^\top)=\sum_{i_1<\dots<i_k} \det(A_{i_1,\dots,i_k})^2$, where 
$A_{i_1,\dots,i_k}$ is the minor of $A$ with indicated rows.
Since the 
$\De_k$ are symmetric we have $\De_k=\tilde\De_k\o \si^n$ for unique 
polynomials $\tilde\De_k$ and similarly we shall use $\tilde B$.

\proclaim{\nmb.{3.2}. Theorem } (Sylvester's version of Sturm's 
theorem, see \cit!{14}, \cit!{10}) 
The roots of $P$ are all 
real if and only if the symmetric $(n\x n)$ matrix $\tilde B(P)$ is 
positive semidefinite; then 
$\tilde\Delta_k(P):=\tilde\De_k(a_1,\dots,a_n)\geq 0$ for $1\le k\le n$. 
The rank of $B$ equals the number 
of distinct roots of $P$ and its signature equals the number of 
distinct real roots. 
\endproclaim 

\proclaim{\nmb.{3.3}. Proposition }
Let now $P$ be a smooth curve of polynomials 
$$ 
P(t)(x)=x^n-a_1(t)x^{n-1}+\dots+(-1)^na_n(t)
$$ 
with all roots real, and distinct for $t=0$.
Then $P$ is smoothly solvable near $0$.

This is also true in the real analytic case and for 
higher dimensional parameters, and in the 
holomorphic case for complex roots.
\endproclaim 

\demo{Proof} 
The derivative $\frac d{dx}P(0)(x)$ does not vanish at any root, 
since they are distinct.
Thus by the implicit function theorem we 
have local smooth solutions $x(t)$ of $P(t,x)=P(t)(x) =0$.
\qed\enddemo 
 
\proclaim{\nmb.{3.4}. Splitting Lemma } 
Let $P_0$ be a polynomial
$$ 
P_0(x)=x^n-a_1x^{n-1}+\dots+(-1)^na_n.
$$ 
If $P_0=P_1\cdot P_2$, 
where $P_1$ and $P_2$ are polynomials 
with no common root. 
Then for $P$ near $P_0$ we have $P=P_1(P)\dot P_2(P)$  
for real analytic mappings of monic polynomials 
$P\mapsto P_1(P)$ and $P\mapsto P_2(P)$, 
defined for $P$ near $P_0$, with the given initial values.
\endproclaim
 
\demo{Proof} 
Let the polynomial $P_0$ be represented as 
the product 
$$ 
P_0=P_1.P_2
=(x^p-b_1x^{p-1}+\dots+(-1)^pb_p)(x^q-c_1x^{q-1}+\dots+(-1)^qc_q),
$$ 
let $x_i$ for $i=1,\dots,n$ be the roots of $P_0$, ordered in such a 
way that for $i=1,\dots,p$ we get the roots of $P_1$, and for 
$i=p+1,\dots,p+q=n$ we get those of $P_2$. 
Then $(a_1,\dots,a_n)=\phi^{p,q}(b_1,\dots,b_p,c_1,\dots,c_q)$
for a polynomial mapping $\phi^{p,q}$ and we get 
$$\align
\si^n&=\phi^{p,q}\o(\sigma^p\times\sigma^q),\\
\det(d\sigma^n)&=\det(d\phi^{p,q}(b,c))\det(d\sigma^p)\det(d\sigma^q).\\
\endalign$$ 
 From \nmb!{3.1} we conclude
$$
\prod_{1\le i<j\le n}\negthickspace(x_i-x_j) 
= \det(d\phi^{p,q}(b,c))
     \prod_{1\le i<j\le p}\negthickspace(x_i-x_j)
     \prod_{p+1\le i<j\le n}\negthickspace(x_i-x_j)
$$
which in turn implies
$$
\det(d\phi^{p,q}(b,c)) = 
\prod_{1\le i\le p<j\le n}\negthickspace(x_i-x_j)\ne 0
$$ 
so that $\phi^{p,q}$ is a real analytic diffeomorphism near $(b,c)$.
\qed  \enddemo 

\subhead\nmb.{3.5}  \endsubhead 
For a continuous function $f$ defined near 0 in $\Bbb R$ let the 
\idx{\it multiplicity} or \idx{\it order of flatness} $m(f)$ at 0 be 
the supremum of all integer $p$ such that
$f(t)=t^pg(t)$ near $0$ for a continuous function $g$.
If f is $C^n$ and 
$m(f)< n$ then $f(t)=t^{m(f)}g(t)$ where now $g$ is $C^{n-m(f)}$
and $g(0)\ne 0$.
If $f$ is a continuous function on the space of polynomials, then for a 
fixed continuous curve $P$ of polynomials we will denote by $m(f)$ the 
multiplicity at 0 of $t\mapsto f(P(t))$.

The splitting lemma \nmb!{3.4} shows that for the problem of smooth 
solvability it is enough to assume that all roots of $P(0)$ are 
equal.
 
\proclaim{Proposition } 
Suppose that the smooth curve of polynomials
$$ 
P(t)(x)=x^n+a_2(t)x^{n-2}-\dots+(-1)^na_n(t)
$$ 
is smoothly solvable with smooth roots $t\mapsto x_i(t)$, 
and that all roots of $P(0)$ are equal.
Then for $(k=2,\dots,n)$
$$\align
m(\tilde\De_k)&\geq k(k-1)\min_{1\le i\le n}m(x_i)\\
m(a_k)&\geq k\min_{1\le i\le n}m(x_i)\\
\endalign$$ 
This result also holds in the real analytic case and in the smooth 
case.
\endproclaim 

\demo{Proof} This follows from 
\nmb!{3.1}.(2) for $\Delta_k$, and from 
$a_k(t)=\si_k(x_1(t),\dots,x_n(t))$. 
\qed\enddemo
 
\proclaim{\nmb.{3.6}. Lemma } Let $P$ be a polynomial of degree $n$ 
with all roots real.
If $a_1=a_2=0$ then all roots of 
$P$ are equal to zero.  
\endproclaim 

\demo{Proof} From \nmb!{3.1}.(1) 
we have $\sum x_i^2=s_2(x)=\sigma_1^2(x)-2\sigma_2(x)=a_1^2-2a_2=0$.
\qed\enddemo 

\proclaim{\nmb.{3.7}. Multiplicity lemma} 
Consider a smooth curve of polynomials
$$ 
P(t)(x)=x^n+a_2(t)x^{n-2}-\dots+(-1)^na_n(t)
$$ 
with all roots real.
Then, for integers $r$, the following conditions are equivalent: 
\roster 
\item  $m(a_k)\geq kr$  for all $2\le k\le n$.
\item $m(\tilde\De_k)\geq k(k-1)r$ for all $2\le k\le n$. 
\item  $m(a_2)\geq 2r$. 
\endroster 
\endproclaim 

\demo{Proof} We only have to treat $r>0$.

(1) implies (2): 
 From \nmb!{3.1}.(1) we have $m(\tilde s_k)\geq rk$, and from the
definition of $\tilde\Delta_k=\det (\tilde B_k)$ we get (2). 

(2) implies (3) since $\tilde\De_2=-2na_2$.

(3) implies (1): 
 From $a_2(0)=0$ and lemma \nmb!{3.6} it follows that 
all roots of the polynomial $P(0)$ are equal to zero and, then, 
$a_3(0)=\dots =a_n(0)=0$.
Therefore, $m(a_3),\dots, m(a_n)\geq 1$. 
Under these conditions, we have 
$a_2(t)=t^{2r}a_{2,2r}(t)$ and $a_k(t)=t^{m_k}a_{k,m_k}(t)$ for 
$k=3,\dots,n$, where the  
$m_k$ are positive integers and $a_{2,2r},a_{3,m_3},\dots,a_{n,m_n}$ 
are smooth functions, and where we may assume that either 
$m_k=m(a_k)<\infty$ or $m_k\ge kr$. 

Suppose now indirectly that for some $k> 2$ we have $m_k=m(a_k)<kr$.
Then we put 
$$ 
m:=\operatorname{min}(r,\frac{m_3}{3},\dots,\frac{m_n}{n}) < r.
$$ 
We consider the following continuous curve of polynomials for 
$t\ge0$:
$$\multline 
\bar P_m(t)(x):=x^n+a_{2,2r}(t)t^{2r-2m}x^{n-2}\\ 
-a_{3,m_3}(t)t^{m_3-3m}x^{n-3}+\dots+ 
(-1)^na_{n,m_n}(t)t^{m_n-nm}. 
\endmultline $$ 
If $x_1,\dots,x_n$ are the real roots of $P(t)$ then 
$t^{-m}x_1,\dots,t^{-m}x_n$ are the roots of $\bar P_m(t)$, for $t>0$.
So for $t>0$, $\bar P_m(t)$ is a family of polynomials with all  
roots real.
Since by theorem \nmb!{3.2} the set of polynomials with 
all roots real is closed, 
$\bar P_m(0)$ is also a polynomial with all roots real. 
 
By lemma \nmb!{3.6} all roots of the polynomial $\bar P_m(0)$ are equal 
to zero, and for those $k$ with $m_k=km$ we have
$a_{k,m_k}(0)=0$ and, therefore, $m(a_k)>m_k$, a contradiction. 
\qed \enddemo 
 
\subhead\nmb.{3.8}. Algorithm \endsubhead
Consider a smooth curve of polynomials
$$ 
P(t)(x)=x^n-a_1(t)x^{n-1}+a_2(t)x^{n-2}-\dots+(-1)^na_n(t)
$$ 
with all roots real.
The algorithm has the following steps:
\roster
\item If all roots of $P(0)$ are pairwise different, $P$ is smoothly 
       solvable for $t$ near 0 by \nmb!{3.3}.
\item If there are distinct roots at $t=0$ we put them into two 
       subsets which splits $P(t)=P_1(t).P_2(t)$ by the splitting lemma 
       \nmb!{3.4}. 
       We then feed $P_i(t)$ (which have lower degree) 
       into the algorithm.
\item All roots of $P(0)$ are equal.
       We first reduce $P(t)$ to the case $a_1(t)=0$ by replacing the 
       variable $x$ by $y=x-a_1(t)/n$.
       Then all roots are equal to 0 so $m(a_2)>0$. 
\item"(3a)"
       If $m(a_2)$ is finite then it is even since 
       $\tilde\De_2=-2na_2\ge0$, $m(a_2)=2r$ and by the multiplicity 
       lemma \nmb!{3.7} 
       $a_i(t)=a_{i,ir}(t)t^{ir}$ $(i=2,\dots,n)$ for smooth $a_{i,ir}$.
       Consider the following smooth curve of polynomials 
$$ 
P_r(t)(x)=x^n+a_{2,2r}(t)x^{n-2}-a_{3,3r}(t)x^{n-3}\dots 
+(-1)^na_{n,nr}(t).
$$ 
       If $P_r(t)$ is smoothly solvable and $x_k(t)$ are its smooth 
       roots, then $x_k(t)t^r$ are the roots of $P(t)$ and the 
       original curve $P$ is smoothly solvable too. 
       Since $a_{2,2m}(0)\ne0$, not all roots of $P_r(0)$ are equal 
       and we may feed $P_r$ into step 2 of the algorithm.
\item"(3b)"
       If $m(a_2)$ is infinite and $a_2=0$, then all roots are 0 by 
       \nmb!{3.6} and thus the polynomial is solvable.
\item"(3c)"
       But if $m(a_2)$ is infinite and  $a_2\ne 0$, then by the 
       multiplicity lemma \nmb!{3.7} all $m(a_i)$ for $2\le i\le n$ 
       are infinite. 
       In this case we keep $P(t)$
       as factor of the original curve of polynomials with all 
       coefficients infinitely flat at $t=0$ after forcing $a_1=0$. 
       This means that all roots of $P(t)$ meet of infinite order of 
       flatness (see \nmb!{3.5}) at $t=0$ for any choice of the 
       roots. 
       This can be seen as follows: If $x(t)$ is any root of 
       $P(t)$ then $y(t):=x(t)/t^r$ is a root of $P_r(t)$, hence by 
       \nmb!{4.1} bounded, so $x(t)=t^{r-1}.ty(t)$ and 
       $t\mapsto ty(t)$ is continuous at $t=0$.
\endroster

This algorithm produces a splitting of the original polynomial 
$$P(t)=P^{(\infty)}(t)P^{(s)}(t)$$
where $P^{(\infty)}$ has the property that each root meets another 
one of infinite order at $t=0$,
and where $P^{(s)}(t)$ is smoothly solvable, and no two roots meet 
of infinite order at $t=0$, if they are not equal.
Any two choices of 
smooth roots of $P^{(s)}$ differ by a permutation.

The factor $P^{(\infty)}$ may or may not be smoothly solvable.
For a flat function $f\ge0$ consider:
$$x^4 - (f(t)+t^2)x^2 +t^2f(t) = (x^2-f(t)).(x-t)(x+t).$$
Here the algorithm produces this factorization.
For $f(t)=g(t)^2$ the polynomial is smoothly solvable.
For the smooth function $f$ from \nmb!{2.4} it is not smoothly 
solvable.

\head\totoc\nmb0{4}. Choosing global roots of polynomials 
differentiably 
\endhead

\proclaim{\nmb.{4.1}. Lemma}
For a polynomial
$$ 
P(x)=x^n-a_1(P)x^{n-1}+\dots+(-1)^na_n(P)
$$ 
with all roots real, i.e\. 
$\tilde\De_k(P)=\tilde\De_k(a_1,\dots,a_n)\ge0$ 
for $1\le k\le n$, let 
$$x_1(P)\le x_2(P)\le \dots \le x_n(P)$$
be the roots, increasingly ordered. 

Then all $x_i:\si^n(\Bbb R^n)\to \Bbb R$ are continuous.
\endproclaim

\demo{Proof}
We show first that $x_1$ is continuous. 
Let $P_0\in \si^n(\Bbb R^n)$ be arbitrary. 
We have to show that for every $\ep >0$
there exists some $\de>0$ such that for all $|P-P_0|<\de$
there is a root $x(P)$ of $P$ with $x(P)<x_1(P_0)+\ep$ and for all roots
$x(P)$ of $P$ we have $x(P)>x_1(P_0)-\ep$.
Without loss of
generality we may assume that $x_1(P_0)=0$.

We use induction on the degree $n$ of $P$. 
By the splitting lemma 
\nmb!{3.4} for the $C^0$-case we may factorize $P$ as
$P_1(P)\cdot P_2(P)$, where $P_1(P_0)$ has all roots equal to $x_1=0$
and $P_2(P_0)$ has all roots greater than $0$ and both polynomials
have coefficients which depend real analytically on $P$.
The degree of $P_2(P)$ is now smaller than $n$, so by induction the 
roots of $P_2(P)$ are continuous and thus larger than 
$x_1(P_0)-\ep$ for $P$ near $P_0$. 

Since $0$ was the smallest root of $P_0$ we have to show
that for all $\ep>0$ there exists a $\de>0$ such that for $|P-P_0|< 
\de$ any
root $x$ of $P_1(P)$ satisfies $|x|<\ep$.
Suppose there is a root $x$ with $|x|\geq \ep$.
Then we get as follows a contradiction, where $n_1$ is the degree of 
$P_1$. 
From 
$$
-x^{n_1} = \sum_{k=1}^{n_1} (-1)^ka_k(P_1) x^{{n_1}-k} 
$$
we have
$$
\ep \leq |x| = \Bigl|\sum_{k=1}^{n_1} (-1)^ka_k(P_1) x^{1-k} \Bigr| 
     \leq \sum_{k=1}^{n_1} |a_k(P_1)|\, |x|^{1-k}  
     < \sum_{k=1}^{n_1} \frac{\ep^k}{{n_1}}\, \ep^{1-k} =\ep , 
$$
provided that ${n_1}|a_k(P_1)|<\ep^k$, which is true for $P_1$ near 
$P_0$, since $a_k(P_0)=0$.
Thus $x_1$ is continuous.

Now we factorize $P=(x-x_1(P))\cdot P_2(P)$, where $P_2(P)$ has roots
$x_2(P)\leq \dots \leq x_n(P)$.
By Horner's algorithm 
($a_n=b_{n-1}x_1$, $a_{n-1}=b_{n-1}+b_{n-2}x_1$,\dots,
$a_2=b_2+b_1x_1$, $a_1=b_1+x_1$)
the coefficients $b_k$ of $P_2(P)$ are again continuous
and so we may proceed by induction on the degree of $P$.
Thus the claim is proved.
\qed\enddemo

\proclaim{\nmb.{4.2}. Theorem} 
Consider a smooth curve of polynomials
$$ 
P(t)(x)=x^n+a_2(t)x^{n-2}-\dots+(-1)^na_n(t)
$$ 
with all roots real, for $t\in \Bbb R$.
Let one of the two following 
equivalent conditions be satisfied:
\roster
\item If two of the increasingly ordered continuous roots meet of 
       infinite order somewhere then they are equal everywhere.   
\item Let $k$ be maximal with the property that $\tilde\De_k(P)$ does not 
       vanish identically for all $t$. 
       Then $\tilde\De_k(P)$ vanishes 
       nowhere of infinite order. 
\endroster

Then the roots of $P$ can be chosen smoothly, and any two choices 
differ by a permutation of the roots.
\endproclaim 

\demo{Proof. The local situation}
We claim that for any $t_0$, without loss $t_0=0$, 
the following conditions are equivalent:
\roster
\item If two of the increasingly ordered continuous roots meet of 
       infinite order at $t=0$ then their germs at $t=0$ are equal.   
\item Let $k$ be maximal with the property that the germ at $t=0$ of 
       $\tilde\De_k(P)$ is not 0. 
       Then $\tilde\De_k(P)$ is not infinitely flat at $t=0$. 
\item The algorithm \nmb!{3.8} never leads to step~(3c).
\endroster
\noindent\therosteritem3 $\Longrightarrow$ \therosteritem1.
Suppose indirectly that two nonequal of the increasingly ordered continuous 
roots meet of infinite order at $t=0$.
Then in each application of 
step~\therosteritem2 these two roots stay with the same factor. 
After any application of step~(3a) these two roots lead to 
nonequal roots of the modified polynomial which still meet of 
infinite order at $t=0$.
They never end up in a facter leading to 
step~(3b) or step~\therosteritem1.
So they end up in a factor leading 
to step~(3c).

\noindent\therosteritem1 $\Longrightarrow$ \therosteritem2.
Let $x_1(t)\le\dots\le x_n(t)$ be the continuous roots of $P(t)$.
 From \nmb!{3.1},~\thetag2 we have
$$ 
\tilde\Delta_k(P(t))
     =\sum_{i_1<i_2<\dots<i_k}(x_{i_1}-x_{i_2})^2\dots(x_{i_1}-x_{i_n})^2 
     \dots(x_{i_{k-1}}-x_{i_k})^2.\tag4
$$ 
The germ of $\tilde\Delta_k(P)$ is not 0, so the germ of one 
summand is not 0.
If $\tilde\Delta_k(P)$ were infinitely flat at $t=0$, then each 
summand is infinitely flat, so there are two roots among the $x_i$ 
which meet of infinite order, thus by assumption their germs are 
equal, so this summand vanishes. 

\noindent\therosteritem2 $\Longrightarrow$ \therosteritem3.
Since the leading $\tilde\De_k(P)$ vanishes only of finite order at 
zero, $P$ has exactly $k$ different roots off 0.
Suppose indirectly that the algorithm \nmb!{3.8} leads to step~(3c), 
then $P=P^{(\infty)}P^{(s)}$ for a nontrivial polynomial 
$P^{(\infty)}$.
Let $x_1(t)\le\dots\le x_p(t)$ be the roots of 
$P^{(\infty)}(t)$ and $x_{p+1}(t)\le\dots\le x_n(t)$ those of 
$P^{(s)}$.
We know that each $x_i$ meets some $x_j$ of infinite 
order and does not meet any $x_l$ of infinite 
order, for $i,j\le p<l$.
Let $k^{(\infty)}>2$ and $k^{(s)}$ be the 
number of generically different roots of $P^{(\infty)}$ and 
$P^{(s)}$, respectively.
Then $k=k^{(\infty)}+k^{(s)}$, and an 
inspection of the formula for $\tilde\De_k(P)$ above leads to the 
fact that it must vanish of infinite order at 0, since the only 
non-vanishing summands involve exactly $k^{(\infty)}$ many 
generically  different roots from $P^{(\infty)}$.

Let $y_1(t)\le\dots\le y_k(t)$ be the 
 From \nmb!{3.1},~thetag2 we 

{\it The global situation.}
 From the first part of the proof we see that the algorithm \nmb!{3.8} 
allows to choose the roots 
smoothly in a neighborhood of each point $t\in \Bbb R$, and that any 
two choices differ by a (constant) permutation of the roots.
Thus we 
may glue the local solutions to a global solution. 
\qed \enddemo 

\proclaim{\nmb.{4.3}. Theorem}
Consider a curve of polynomials
$$ 
P(t)(x)=x^n-a_1(t)x^{n-1}+\dots+(-1)^na_n(t), \quad t\in \Bbb R,
$$ 
with all roots real, where all $a_i$ are of class $C^n$.
Then there is a differentiable curve 
$x=(x_1,\dots,x_n):\Bbb R\to\Bbb R^n$ whose coefficients parameterize 
the roots. 
\endproclaim

That this result cannot be improved to $C^2$-roots is shown already 
in \nmb!{2.4}, and not to $C^1$ for $n\ge 3$ is shown in \nmb!{2.5}.

\demo{Proof}
First we note that the multiplicity lemma \nmb!{3.7} remains true in 
the $C^n$-case for $r=1$ in the following sense, with the same proof:
\newline
{\sl If $a_1=0$ then the following two conditions are equivalent
\roster
\item $a_k(t)=t^k a_{k,k}(t)$ for a continuous function $a_{k,k}$, 
       for $2\le k\le n$.
\item $a_2(t)=t^2 a_{2,2}(t)$ for a continuous function $a_{2,2}$.
\endroster}

In order to prove the theorem itself we follow one step of the 
algorithm.
First we replace $x$ by $x+\frac1na_1(t)$, or assume 
without loss that $a_1=0$.
Then we choose a fixed $t$, say $t=0$. 

If $a_2(0)=0$ then it vanishes of second order at 0: if it vanishes only 
of first order then $\tilde\De_2(P(t))=-2na_2(t)$ would change sign 
at $t=0$, contrary to the assumption that all roots of $P(t)$ are 
real, by \nmb!{3.2}.
Thus $a_2(t)=t^2a_{2,2}(t)$, so by the variant 
of the multiplicity lemma \nmb!{3.7} described above we have 
$a_k(t)=t^k a_{k,k}(t)$ for continuous functions $a_{k,k}$, for 
$2\le k\le n$. 
We consider the following continuous curve of polynomials 
$$ 
P_1(t)(x)=x^n+a_{2,2}(t)x^{n-2}-a_{3,3}(t)x^{n-3}\dots 
+(-1)^na_{n,n}(t).
$$ 
with continuous roots $z_1(t)\le\dots\le z_n(t)$, by \nmb!{4.1}. 
Then $x_k(t)=z_k(t)t$ are differentiable at $0$, and are all 
roots of $P$, but note that $x_k(t)=y_k(t)$ for $t\ge0$, but 
$x_k(t)=y_{n-k}(t)$ for $t\le 0$, where $y_1(t)\le\dots\le y_n(t)$ 
are the ordered roots of $P(t)$.
This gives us one choice of differentiable roots near $t=0$.
Any choice is then given by this choice and applying afterwards
any permutation of the set $\{1,\dots,n\}$ keeping invariant the 
function $k\mapsto z_k(0)$. 

If $a_2(0)\ne0$ then by the splitting lemma \nmb!{3.4} for the $C^n$-case 
we may factor $P(t)=P_1(t)\dots P_k(t)$ where 
the $P_i(t)$ have again $C^n$-coefficients and 
where each $P_i(0)$ has all roots equal to $c_i$, and where the $c_i$ 
are distinct.
By the arguments above the roots of each $P_i$ can be 
arranged differentiably, thus $P$ has differentiable roots $y_k(t)$. 

But note that we have to apply a permutation on one side of 0 to 
the original roots, in the following case:
Two roots $x_k$ and $x_l$ meet at zero with 
$x_k(t)-x_l(t)=tc_{kl}(t)$ with $c_{kl}(0)\ne0$ 
(we say that they meet slowly).
We may apply to this choice an arbitrary 
permutation of any two roots which meet with $c_{kl}(0)=0$ (i.e\. 
at least of second order), and we get thus any differentiable
choice near $t=0$.

Now we show that we choose the roots differentiable on the whole 
domain $\Bbb R$.
We start with with the ordered continuous roots
$y_1(t)\le\dots\le y_n(t)$. 
Then we put 
$$x_k(t)=y_{\si(t)(k)}(t)$$
where the permutation $\si(t)$ is given by 
$$\si(t)=(1,2)^{\ep_{1,2}(t)}\dots(1,n)^{\ep_{1,n}(t)}(2,3)^{\ep_{2,3}(t)}
     \dots(n-1,n)^{\ep_{n-1,n}(t)}$$
and where $\ep_{i,j}(t)\in\{0,1\}$ will be specified as follows:
On the closed set $S_{i,j}$ of all $t$ where $y_i(t)$ and $y_j(t)$ 
meet of order at least 2 any choice is good.
The complement of 
$S_{i,j}$ is an at most countable union of open intervals, and in each 
interval we choose a point where we put $\ep_{i,j}=0$.
Going right 
(and left) from this point we change $\ep_{i,j}$ in each point where 
$y_i$ and $y_j$ meet slowly.
These points accumulate only in 
$S_{i,j}$. 
\qed\enddemo

\head\totoc\nmb0{5}. The real analytic case \endhead

\proclaim{\nmb.{5.1}. Theorem}
Let $P$ be a real analytic curve of polynomials 
$$ 
P(t)(x)=x^n-a_1(t)x^{n-1}+\dots+(-1)^na_n(t),\quad t\in \Bbb R,
$$ 
with all roots real. 

Then $P$ is real analytically solvable, globally on $\Bbb R$.
All solutions differ by permutations.
\endproclaim

By a real analytic curve of polynomials we mean that all $a_i(t)$ are 
real analytic in $t$ (but see also \cit!{8}), 
and real analytically solvable means that we 
may find $x_i(t)$ for $i=1,\dots,n$ which are real 
analytic in $t$ and are roots of $P(t)$ for all $t$. 
The local existence part of this theorem is due to Rellich \cit!{11}, 
Hilfssatz 2, his proof uses Puiseux-expansions.
Our proof is 
different and more elementary.

\demo{Proof}
We first show that $P$ is  locally real analytically solvable, near 
each point $t_0\in \Bbb R$.
It suffices to consider $t_0=0$. 
Using the transformation in the introduction we first assume that 
$a_1(t)=0$ for all $t$.
We use induction on the degree $n$.
If $n=1$ the theorem holds.
For $n>1$ we consider several cases:

The case $a_2(0)\ne0$.
Here not all roots of $P(0)$ are 
equal and zero, so by the splitting lemma \nmb!{3.4} we may factor 
$P(t)=P_1(t).P_2(t)$ for real analytic curves of polynomials of 
positive degree, which have both all roots real, and we have reduced 
the problem to lower degree.

The case $a_2(0)=0$. 
If $a_2(t)=0$ for all $t$, then by \nmb!{3.6} 
all roots of $P(t)$ are 0, and we are done. 
Otherwise $1\le m(a_2)<\infty$ for the multiplicity of $a_2$ at 0, 
and by \nmb!{3.6} all roots of $P(0)$ are 0. 
If $m(a_2)>0$ is odd, then $\tilde\De_2(P)(t)=-2na_2(t)$ changes sign 
at $t=0$, so by \nmb!{3.2} not all roots of $P(t)$ are real for $t$ 
on one side of 0.
This contradicts the assumption, so $m(a_2)=2r$ 
is even.
Then by the multiplicity lemma \nmb!{3.7} we have 
$a_i(t)=a_{i,ir}(t)t^{ir}$ $(i=2,\dots,n)$ for real analytic 
$a_{i,ir}$, and we may consider the following real analytic curve of 
polynomials
$$ 
P_r(t)(x)=x^n+a_{2,2r}(t)x^{n-2}-a_{3,3r}(t)x^{n-3}\dots 
+(-1)^na_{n,nr}(t),
$$
with all roots real.
If $P_r(t)$ is real analytically solvable and 
$x_k(t)$ are its real analytic roots, then $x_k(t)t^r$ are the roots 
of $P(t)$ and the original curve $P$ is real analytically solvable 
too.
Now $a_{2,2r}(0)\ne 0$ and we are done by the case above.

{\sl Claim.}
Let $x=(x_1,\dots,x_n):I\to \Bbb R^n$
be a real analytic curve of roots of $P$ on an open interval 
$I\subset \Bbb R$.
Then any real analytic curve of roots of $P$ on 
$I$ is of the form $\al\o x$ for some permutation $\al$.

Let $y:I\to \Bbb R^n$ be another real analytic curve of roots of $P$.
Let $t_k\to t_0$ be a convergent sequence of distinct points in $I$.
Then $y(t_k)=\al_k(x(t_k))=(x_{\al_k1},\dots,x_{\al_kn})$ 
for permutations $\al_k$.
By choosing a subsequence we may assume 
that all $\al_k$ are the same permutation $\al$.
But then the real analytic curves $y$ and $\al\o x$ coincide on a 
converging sequence, so they coincide on $I$ and the claim follows.

Now from the local smooth solvability above and the uniqueness of 
smooth solutions up to permutations we can glue a global smooth 
solution on the whole of $\Bbb R$.
\qed\enddemo

\subhead\nmb.{5.2}. Remarks and examples \endsubhead
The uniqueness statement of theorem \nmb!{5.1} is wrong in the smooth 
case, as is shown by the following  
example: $x^2=f(t)^2$ where $f$ is smooth.
In each point $t$ where 
$f$ is infinitely flat one can change sign in the solution 
$x(t)=\pm f(t)$.
No sign change can be absorbed in a permutation 
(constant in $t$).
If there are infinitely many points of flatness 
for $f$ we get uncountably many smooth solutions.

Theorem \nmb!{5.1} reminds of the curve lifting property of covering 
mappings.
But unfortunately one cannot lift real analytic homotopies, 
as the following example shows.
This example also shows that polynomials which are real analytically 
parameterized by higher dimensional variables are not real 
analytically solvable.

Consider the 2-parameter family $x^2=t_1^2+t_2^2$. 
The two continuous solutions are 
$x(t)=\pm|t|$, but for none of them 
$t_1\mapsto x(t_1,0)$ is differentiable at 0. 

There remains the question whether for a real analytic submanifold 
of the space of polynomials with all roots real one can choose 
the roots real analytically along this manifold.
The following 
example shows that this is not the case:

Consider 
$$P(t_1,t_2)(x)=(x^2-(t_1^2+t_2^2))\,(x-(t_1-a_1))\,(x-(t_2-a_2)),$$
which is not real analytically solvable, see above.
For $a_1\ne a_2$ 
the coefficients describe a real analytic embedding for $(t_1,t_2)$ 
near 0.

\head\totoc \nmb0{6}. Choosing roots of complex polynomials 
\endhead 

\subhead\nmb.{6.1}  \endsubhead
In this section we consider the problem of finding smooth curves of 
complex roots for smooth curves $t\mapsto P(t)$ of polynomials 
$$ 
P(t)(z)=z^n-a_1(t)z^{n-1}+\dots+(-1)^na_n(t) 
$$ 
with complex valued coefficients $a_1(t),\dots,a_n(t)$. 
We shall also discuss the real analytic and holomorphic cases. 
The definition of the Bezoutiant $B$, its 
principal minors $\Delta_k$, and formula \nmb!{3.1}.(2) are still 
valid.
Note that now there are no restrictions on the coefficients. 
In this section the parameter may be real and $P(t)$ may be smooth or 
real analytic in $t$, or the parameter $t$ may be complex and $P(t)$ 
may be holomorphic in $t$.
 
\subhead\nmb.{6.2}. The case $n=2$ \endsubhead
As in the real 
case the problem reduces to the following one: 
Let $f$ be a smooth complex valued function, defined near $0$ in 
$\Bbb R$, such that $f(0)=0$.
We look for a smooth function 
$g:(\Bbb R,0)\to \Bbb C$ such that $f=g^2$.
If $m(f)$ is finite and even, we have $f(t)=t^{m(f)}h(t)$ with 
$h(0)\ne 0$, and $g(t):=t^{m(f)/2}\sqrt{h(t)}$ is a local solution.
If $m(f)$ is finite and odd there is no solution $g$, also not in the 
real analytic and holomorphic cases.
If $f(t)$ is flat at $t=0$, then one has no definite answer, and the 
example \nmb!{2.4} is still not smoothly solvable.

\subhead\nmb.{6.3}. The general case \endsubhead
Proposition \nmb!{3.3} and the splitting lemma \nmb!{3.4} are true in the 
complex case.
Proposition \nmb!{3.5} is true also because it follows 
from \nmb!{3.1}.(2).
Of course lemma \nmb!{3.6} is not true now and 
the multiplicity lemma \nmb!{3.7} only partially holds:

\proclaim{\nmb.{6.4}. Multiplicity Lemma} 
Consider the smooth (real analytic, holomorphic) 
curve of complex polynomials
$$ 
P(t)(z)=z^n+a_2(t)z^{n-2}-\dots+(-1)^na_n(t).
$$ 
Then, for integers $r$, the following conditions are equivalent: 
\roster 
\item  $m(a_k)\geq kr$  for all $2\le k\le n$.
\item $m(\tilde\De_k)\geq k(k-1)r$ for all $2\le k\le n$. 
\endroster 
\endproclaim 

\demo{Proof} (2) implies (1): 
Since $\tilde\De_2=na_2$ we have $s_2(0)=-2a_2(0)=0$.  
 From $\tilde\De_3(0)=-s_3(0)^2$ we then get $s_3(0)=0$, and so on  
we obtain $s_4(0)=\dots=s_n(0)=0$.
Then by \nmb!{3.1}.(1) 
$a_i(0)=0$ for $i=3,\dots,n$.
The rest of the proof coincides with 
the one of the multiplicity lemma \nmb!{3.7}. 
\qed\enddemo 

\subhead\nmb.{6.5}. Algorithm \endsubhead
Consider a smooth (real analytic, holomorphic) curve of polynomials
$$ 
P(t)(z)=z^n-a_1(t)z^{n-1}+a_2(t)z^{n-2}-\dots+(-1)^na_n(t)
$$ 
with complex coefficients.
The algorithm has the following steps:
\roster
\item If all roots of $P(0)$ are pairwise different, $P$ is smoothly 
       (real analytically, holomorphically)  
       solvable for $t$ near 0 by \nmb!{3.3}.
\item If there are distinct roots at $t=0$ we put them into two 
       subsets which factors $P(t)=P_1(t).P_2(t)$ by the splitting lemma 
       \nmb!{3.4}. 
       We then feed $P_i(t)$ (which have lower degree) 
       into the algorithm.
\item All roots of $P(0)$ are equal.
       We first reduce $P(t)$ to the case $a_1(t)=0$ by replacing the 
       variable $x$ by $y=x-a_1(t)/n$.
       Then all roots are equal to 0  
       so $a_i(0)=0$ for all $i$. 
 
       If there does not exist an integer $r>0$ with $m(a_i)\ge ir$ 
       for $2\le i\le n$, then by \nmb!{3.5} the polynomial is 
       {\it not} 
       smoothly (real analytically, holomorphically) solvable, by 
       proposition \nmb!{3.5}: We store the polynomial as an output 
       of the procedure, as a factor of $P^{(n)}$ below.

       If there exists an integer $r>0$ with $m(a_i)\ge ir$ for 
       $2\le i\le n$, let   
       $a_i(t)=a_{i,ir}(t)t^{ir}$ $(i=2,\dots,n)$ for smooth (real 
       analytic, holomorphic) $a_{i,ir}$. 
       Consider the following smooth (real analytic, holomorphic) 
       curve of polynomials  
$$ 
P_r(t)(x)=x^n+a_{2,2r}(t)x^{n-2}-a_{3,3r}(t)x^{n-3}\dots 
+(-1)^na_{n,nr}(t). 
$$ 
       If $P_r(t)$ is smoothly (real analytically, holomorphically) 
       solvable and $x_k(t)$ are its smooth (real analytic, holomorphic) 
       roots, then $x_k(t)t^r$ are the roots of $P(t)$ and the 
       original curve $P$ is smoothly (real analytically, 
       holomorphically) solvable too.  

       If for one coefficient we have $m(a_i)=ir$ then $P_r(0)$ has a 
       coefficient which does not vanish, so not all roots of 
       $P_r(0)$ are equal, and we may feed $P_r$ into step 
       \therosteritem2.

       If all coefficients of $P_r(0)$ are zero, we feed $P_r$ again 
       into step \therosteritem3.

       In the smooth case all $m(a_i)$ may be infinite; In this case 
       we store the polynomial as a factor of 
       $P^{(\infty)}$ below.
\endroster

In the holomorphic and real analytic cases the algorithm provides a 
splitting of the polynomial $P(t)=P^{(n)}(t)P^{(s)}(t)$ into 
holomorphic and real analytic curves, where $P^{(s)}(t)$ is solvable, 
and where $P^{(n)}(t)$ is not solvable.
But it may contain solvable 
roots, as is seen by simple examples.

In the smooth case the algorithm provides a splitting near $t=0$
$$
P(t)=P^{(\infty)}(t)P^{(n)}(t)P^{(s)}(t)
$$
into smooth curves of 
polynomials, where: 
$P^{(\infty)}$ has the property that each root meets another 
one of infinite order at $t=0$;
and where $P^{(s)}(t)$ is smoothly solvable, and no two roots meet 
of infinite order at $t=0$;
$P^{(n)}$ is not smoothly solvable, with the same property as above. 

\subhead\nmb.{6.6}. Remarks \endsubhead
If $P(t)$ is a polynomials whose coefficients are meromorphic 
functions of a complex variable $t$, there is a well developed 
theory of the roots of $P(t)(x)=0$ as multi-valued meromorphic 
functions, given by \idx{\it Puiseux} or \idx{\it Laurent-Puiseux} 
series.
But it is difficult to extract holomorphic information out of 
it, and the algorithm above complements this theory. 
See for example Theorem 3 on page 370 (Anhang, \S 5) of \cit!{1}.
The question of choosing roots continuously has been treated in 
\cit!{2}: one finds sufficient conditions for it.

\head\totoc\nmb0{7}. Choosing eigenvalues and eigenvectors of 
matrices smoothly \endhead 

In this section we consider the following situation:
Let $A(t)=(A_{ij}(t))$ be a smooth (real analytic, holomorphic) curve 
of real (complex) $(n\x n)$-matrices or operators, 
depending on a real (complex) parameter $t$ near 0.
What can we say about the eigenvalues and eigenfunctions of $A(t)$?

Let us first recall some known results.
These have some 
difficulty with the interpretation of the eigenprojections and the 
eigennilpotents at branch 
points of the eigenvalues, see \cit!{7}, II,~1.11.

\proclaim{\nmb.{7.1}.  Result} {\rm (\cit!{7}, II,~1.8)}
Let $\Bbb C\ni t\mapsto A(t)$ 
be a holomorphic curve.
Then all eigenvalues, all 
eigenprojections and all eigennilpotents are holomorphic with at most 
algebraic singularities at discrete points.
\endproclaim

\proclaim{\nmb.{7.2}. Result} {\rm (\cit!{11}, Satz 1)}
Let $t\mapsto A(t)$ be a real analytic curve of hermitian complex 
matrices.
Let  
$\la$ be a $k$-fold eigenvalue of $A(0)$ with $k$ orthonormal 
eigenvectors $v_i$, 
and suppose that there is no other eigenvalue of $A$ near $\la$.
Then there are $k$ real analytic eigenvalues $\la_i(t)$ through 
$\la$, and $k$ orthonormal real analytic eigenvectors through the 
$v_i$, for $t$ near 0.
\endproclaim

The condition that $A(t)$ is hermitian cannot be 
omitted.
Consider the following example of real semisimple (not 
normal) matrices
$$\gather
A(t) := \pmatrix 2t + t^3 & t \\
                -t       & 0 \endpmatrix, \\
\la_{\pm}(t) = t + \frac{t^2}2 \pm t^2\sqrt{1+\tfrac{t^2}4},\quad
x_\pm(t) = \pmatrix 1 +\frac{t}2 \pm t\sqrt{1+\tfrac{t^2}4} \\ 
                                        -1  \endpmatrix,
\endgather$$
where at $t=0$ we do not get a base of eigenvectors.

\proclaim{\nmb.{7.3}. Result} 
{\rm (Rellich \cit!{13}, see also Kato \cit!{7}, II,~6.8)}
Let $A(t)$ be a $C^1$-curve of symmetric matrices.
Then the 
eigenvalues can be chosen $C^1$ in $t$, on the whole parameter 
interval.
\endproclaim

 For an extension of this result to Hilbert space, under stronger 
assumptions, see \nmb!{7.8}, whose proof will need \nmb!{7.3}.
This result is best possible for the degree of continuous 
differentiability, as is shown by the following example.

\subhead\nmb.{7.4} Example 
\endsubhead
Consider the symmetric matrix 
$$
A(t) = \pmatrix
	a(t) & b(t) \\
	b(t) & -a(t) 
	\endpmatrix
$$
The characteristic polynomial of $A(t)$	is 
$\la^2 - (a(t)^2+b(t)^2)$.
We shall specify the entries $a$ and $b$ as smooth functions in such a
way, that $a(t)^2+b(t)^2$ does not admit a $C^2$-square root.

Assume that $a(t)^2+b(t)^2=c(t)^2$ for a $C^2$-function $c$.
Then we may
compute as follows:
$$
\align
c^2 &= a^2 + b^2 \\
cc' &= aa' + bb' \\
(c')^2 + cc'' &= (a')^2 + aa'' + (b')^2 + bb'' \\
c'' &= \frac1c \Bigl( (a')^2 + aa'' + (b')^2 + bb'' - (c')^2 \Bigr)\\
&= \frac1c \Bigl( (a')^2 + aa'' + (b')^2 + bb''-\frac{(aa'+bb')^2}{c^2}\Bigr) \\
&= \frac{(ab'-ba')^2+a^3a''+b^3b''+ab^2a''+a^2bb''}{c^3}
\endalign
$$
By $c^2 = a^2 + b^2$ we have
$$
\left|\frac{a^3}{c^3}\right|\le 1, \quad 
\left|\frac{b^3}{c^3}\right|\le 1, \quad
\left|\frac{ab^2}{c^3}\right|\le \frac1{\sqrt3}, 
\quad \left|\frac{a^2b}{c^3}\right|\le \frac1{\sqrt3}. 
$$
So for $C^2$-functions $a$, $b$, and continuous $c$ 
all these terms are bounded.  
We will now construct smooth $a$ and $b$ such that
$$
\left(\frac{(ab'-ba')^2}{c^3}\right)^2=  
\frac{(ab'-ba')^4}{(a^2+b^2)^3}  
$$
is unbounded near $t=0$.
This contradicts that $c$ is $C^2$.
For this we choose $a$ and $b$ similar to the function $f$ in 
example \nmb!{2.4} with the same $t_n$ and $h_n$:
$$\align
a(t) :&= \sum_{n=1}^\infty h_n(t-t_n). 
     \left(\frac{2n}{2^n}(t-t_n)+\frac1{4^n}\right),\\
b(t) :&= \sum_{n=1}^\infty h_n(t-t_n). 
     \left(\frac{2n}{2^n}(t-t_n)\right).
\endalign$$
Then $a(t_n)=\frac1{4^n}$, $b(t_n)=0$, $|c(t_n)|=\frac1{4^n}$, 
and $b'(t_n)=\frac{2n}{2^n}$.

\proclaim{\nmb.{7.5}. Result} 
{\rm (Rellich \cit!{12}, see also Kato \cit!{7}, VII,~3.9)}
Let $A(t)$ be a real analytic curve of unbounded self-adjoint operators in a 
Hilbert space with common domain of definition and with compact 
resolvent. 

Then the eigenvalues and the 
eigenvectors can be chosen real analytically in $t$, on the whole 
parameter domain.
\endproclaim

\proclaim{\nmb.{7.6}. Theorem} 
Let $A(t)=(A_{ij}(t))$ be a smooth curve of complex hermitian
$(n\x n)$-matrices, depending  on a real parameter $t\in \Bbb R$, 
acting on a hermitian space $V=\Bbb C^n$, such that no two 
of the continuous eigenvalues meet of infinite order at any 
$t\in \Bbb R$ if they are not equal for all $t$. 

Then all the eigenvalues and all the 
eigenvectors can be chosen smoothly in $t$, on the whole 
parameter domain $\Bbb R$.
\endproclaim

The last condition permits that some eigenvalues agree for all $t$ 
--- we speak of higher `generic multiplicity' in this situation.

\demo{Proof}
The proof will use an algorithm.

Note first that by \nmb!{4.2} the characteristic 
polynomial 
$$\align
P(A(t))(\la) &= \det (A(t)-\la\Bbb I) \tag1\\
&=\la^n-a_1(t)\la^{n-1}+a_2(t)\la^{n-2}-\dots+(-1)^na_n(t)\\
&= \sum_{i=0}^n \operatorname{Trace}(\La^iA(t))\la^{n-i}
\endalign$$
is smoothly solvable, with smooth roots 
$\la_1(t),\dots\la_n(t)$, on the whole parameter interval. 

\demo{Case 1: distinct eigenvalues}
If $A(0)$ has some eigenvalues distinct, then one can reorder them in 
such a way that for $i_0=0<1\le i_1< i_2< \dots < i_k < n=i_{k+1}$ we 
have 
$$
\la_1(0)=\dots=\la_{i_1}(0) < \la_{i_1+1}(0)=\dots=\la_{i_2}(0) < 
     \dots < \la_{i_k+1}(0)=\dots=\la_{n}(0)
$$
For $t$ near 0 we still have
$$
\la_1(t),\dots,\la_{i_1}(t) < \la_{i_1+1}(t),\dots,\la_{i_2}(t) < 
     \dots < \la_{i_k+1}(t),\dots,\la_{n}(t)
$$
Now for $j=1,\dots,k+1$ consider the subspaces 
$$
V^{(j)}_t =\bigoplus_{i=i_{j-1}+1}^{i_j} 
     \{v\in V: (A(t)-\la_i(t))v=0\} 
$$
{\sl Then each $V^{(j)}_t$ runs through a smooth vector 
subbundle of the trivial bundle $(-\ep,\ep)\x V \to (-\ep,\ep)$, which 
admits a smooth framing $e_{i_{j-1}+1}(t),\dots,e_{i_j}(t)$.
We 
have $V=\bigoplus_{j=1}^{k+1} V^{(j)}_t$ for each $t$.
}\newline
In order to prove this statement note that 
$$
V^{(j)}_t = \ker\Bigl(
(A(t)-\la_{i_{j-1}+1}(t))\o\dots\o (A(t)-\la_{i_j}(t))\Bigr)
$$
so $V^{(j)}_t$ is the kernel of a smooth vector bundle 
homomorphism $B(t)$ of constant rank (even of constant dimension of 
the kernel), and thus is a smooth vector 
subbundle.
This together with a smooth frame field can be 
shown as follows: Choose a basis of $V$ such that 
$A(0)$ is diagonal.
Then by the elimination procedure one can 
construct a basis for the kernel of $B(0)$.
For $t$ near 0, the 
elimination procedure (with the same choices) gives then a basis of the 
kernel of $B(t)$; the elements of this basis are then smooth 
in $t$, for $t$ near 0.

 From the last result it follows that it suffices to find smooth 
eigenvectors in each subbundle $V^{(j)}$ separately, 
expanded in the smooth frame field.
But in this frame field the 
vector subbundle looks again like a constant vector space.
So feed 
each of this parts ($A$ restricted to $V^{(j)}$, as matrix with 
respect to the frame field) into case 2 below.
\enddemo

\demo{Case 2: All eigenvalues at 0 are equal}
So suppose that $A(t):V\to V$ is hermitian with all eigenvalues 
at $t=0$ equal to $\frac{a_1(0)}n$, see \thetag1. 

Eigenvectors of $A(t)$ are also eigenvectors of 
$A(t)-\frac{a_1(t)}n\Bbb I$, so we 
may replace $A(t)$ by $A(t)-\frac{a_1(t)}n\Bbb I$ and assume that 
for the characteristic polynomial \thetag1 we have $a_1=0$, 
or assume without loss that  $\la_i(0)=0$ for all $i$ and so
$A(0)=0$. 

If $A(t)=0$ for all $t$ we choose the eigenvectors constant.

Otherwise let $A_{ij}(t)=tA^{(1)}_{ij}(t)$. 
 From \thetag1 we 
see that the characteristic polynomial of the hermitian matrix
$A^{(1)}(t)$ is $P_1(t)$ in the notation of \nmb!{3.8},
thus $m(a_i)\ge i$ for $2\le i\le n$ (which follows from \nmb!{3.5} 
also). 

The eigenvalues of $A^{(1)}(t)$ are the roots of $P_1(t)$, which 
may be chosen in a smooth way, since they again satisfy the condition 
of theorem \nmb!{4.2}.
Note that eigenvectors of 
$A^{(1)}$ are also eigenvectors of $A$.
If the eigenvalues are still all equal, we 
apply the same procedure again, until they are not all equal: we 
arrive at this situation by the assumption of the theorem.
Then we apply case 1.
\enddemo

This algorithm shows that one may choose the eigenvectors $x_i(t)$ of 
$A(t)$ in a smooth way, locally in $t$.
It remains to extend this 
to the whole parameter interval. 

If some eigenvalues coincide locally, then on the whole of $\Bbb R$, 
by the assumption. 
The corresponding eigenspaces then form a smooth vector bundle 
over $\Bbb R$, by case 1, since those eigenvalues, which meet in isolated 
points are different after application of case 2. 

So we we get $V=\bigoplus W^{(j)}_t$ where each $W^{(j)}_t$ is a   
smooth sub vector bundles of $V\x \Bbb R$, whose dimension is the 
generic multiplicity of the corresponding smooth eigenvalue 
function.
It suffices to find global orthonormal smooth frames 
for each of these vector bundles; 
this exists since the vector bundle is smoothly 
trivial, by using parallel transport with respect to a smooth 
Hermitian connection.  
\qed\enddemo

\subhead\nmb.{7.7}. Example \endsubhead (see \cit!{11}, \S 2)
That the last result cannot be improved 
is shown by the following example which rotates a lot:
$$\allowdisplaybreaks\align
x_+(t) :&= \pmatrix \cos\frac1t \\ \sin\frac1t \endpmatrix, \quad
x_-(t) := \pmatrix -\sin\frac1t \\ \cos\frac1t \endpmatrix, \quad
\la_\pm(t) = \pm e^{-\frac1{t^2}},\\
A(t) :&= (x_+(t),x_-(t))
     \pmatrix \la_+(t) & 0 \\
                 0 & \la_-(t) \endpmatrix
      (x_+(t),x_-(t))\i \\
&= e^{-\frac1{t^2}}\pmatrix \cos\frac2t & \sin\frac2t \\
                         \sin\frac2t & -\cos\frac2t 
                                \endpmatrix.
\endalign$$
Here $t\mapsto A(t)$ and $t\mapsto \la_\pm(t)$ are smooth, whereas 
the eigenvectors cannot be chosen continuously.

\proclaim{\nmb.{7.8}. Theorem}
Let $t\mapsto A(t)$ be a smooth curve of unbounded self-adjoint 
operators in a Hilbert space with common domain of definition and with 
compact resolvent.
Then the eigenvalues of $A(t)$ may be arranged 
in such a way that each eigenvalue is $C^1$.

Suppose moreover that no two of the continuously chosen eigenvalues 
meet of infinite order at any $t\in \Bbb R$ if they are not equal.  
Then the eigenvalues and the eigenvectors can be chosen smoothly in 
$t$, on the whole parameter domain.
\endproclaim

\subhead Remarks \endsubhead
That $A(t)$ is a smooth curve of unbouded operators means the 
following:
There is a dense subspace $V$ of the Hilbert space $H$ 
such that $V$ is the domain of definition of each $A(t)$, and such 
that $A(t)^*=A(t)$ with the same domains $V$, 
where the adjoint operator $A(t)^*$ is defined 
by $\langle A(t)u,v \rangle=\langle u,A(t)^*v\rangle$ 
for all $v$ for which the left hand side 
is bounded as function in $u\in H$.
Moreover we require 
that $t\mapsto \langle A(t)u,v\rangle$ 
is smooth for each $u\in V$ and $v\in H$. 
This implies that $t\mapsto A(t)u$ is smooth $\Bbb R\to H$ for each 
$u\in V$ by \cit!{9},~2.3 or \cit!{5},~2.6.2. 

The first part of the proof will show that $t\mapsto A(t)$ smooth 
implies that the resolvent $(A(t)-z)\i$ is smooth in $t$ and $z$ 
jointly, and only this is used later in the proof.

It is well known and in the proof we will show that if for some 
$(t,z)$ the resolvent $(A(t)-z)\i$ is compact then for all 
$t\in \Bbb R$ and $z$ in the resolvent set of $A(t)$.
 
\demo{Proof} 
For each $t$ consider the norm $\|u\|_t^2:=\|u\|^2+\|A(t)u\|^2$ on 
$V$.
Since $A(t)=A(t)^*$ is closed, $(V,\|\quad\|_t)$ is again a 
Hilbert space with inner product 
$\langle u,v\rangle_t:=\langle u,v\rangle+\langle A(t)u,A(t)v\rangle$. 
All these norms are equivalent since 
$(V,\|\quad\|_t+\|\quad\|_s)\to (V,\|\quad\|_t)$ is continuous and 
bijective, so an isomorphism by the open mapping theorem. 
Then $t\mapsto \langle u,v\rangle_t$ is smooth for fixed $u,v\in V$, 
and by the multilinear uniform boundedness principle 
(\cit!{9},~5.17 or \cit!{5},~3.7.4~+~4.1.19) the mapping 
$t\mapsto \langle \quad,\quad\rangle_t$ is smooth into the space of 
bounded bilinear forms on $(V,\|\quad\|_s)$ for each fixed $s$. 
By the exponential law 
(\cit!{9},~3.12 or \cit!{5},~1.4.3)
$(t,u)\mapsto \|u\|^2_t$ is smooth from 
$\Bbb R\x (V,\|\quad\|_s)\to \Bbb R$ for each fixed $s$.
Thus all Hilbert norms $\|\quad\|_t$ are locally uniformly 
equivalent, since  
$\{\|u\|_t:|t|\le K,\|u\|_s\le 1 \}$ is bounded by $L_{K,s}$ in 
$\Bbb R$, so $\|u\|_t\le L_{K,s}\|u\|_s$ for all $|t|\le K$. 
Let us now equip $V$ with one of the equivalent Hilbert norms, say 
$\|\quad\|_0$.
Then each $A(t)$ is a globally defined operator $V\to H$ 
with closed graph and is thus bounded, and by using again the 
(multi)linear uniform boundedness theorem as above we see
that $t\mapsto A(t)$ is smooth $\Bbb R\to L(V,H)$.

If for some $(t,z)\in \Bbb R\x\Bbb C$ the bounded operator 
$A(t)-z:V\to H$ is invertible, then this is true locally and 
$(t,z)\mapsto (A(t)-z)\i:H\to V$ is smooth since inversion is smooth 
on Banach spaces.

Since each $A(t)$ is hermitian the \idx{\it global resolvent set}
$\{(t,z)\in \Bbb R\x\Bbb C: (A(t)-z):V\to H \text{ is invertible}\}$  
is open, contains $\Bbb R\x(\Bbb C\setminus \Bbb R)$, and hence is 
connected. 

Moreover $(A(t)-z)\i:H\to H$ is a compact operator 
for some (equivalently any) $(t,z)$ if and only if the inclusion 
$i:V\to H$ is compact, since $i=(A(t)-z)\i\o(A(t)-z): V\to H\to H$.

Let us fix a parameter $s$. 
We choose a simple smooth curve $\ga$ in the resolvent set of $A(s)$ for 
fixed $s$. 
\roster
\item {\it Claim} For $t$ near $s$, there are $C^1$-functions 
       $t\mapsto\la_i(t):1\le i\le N$ which parametrize all 
       eigenvalues (repeated according to their multiplicity) 
       of $A(t)$ in the interior of $\ga$. 
       If no two of the 
       generically different eigenvalues meet of infinite order they 
       can be chosen smoothly. 
\endroster
By replacing $A(s)$ by $A(s)-z_0$ if 
necessary we may assume that 0 is not an eigenvalue of $A(s)$. 
Since the global resolvent set is open, no eigenvalue 
of $A(t)$ lies on $\ga$ or equals 0, for $t$ near $s$.
Since 
$$
t\mapsto -\frac1{2\pi i}\int_\ga (A(t)-z)\i\;dz =: P(t,\ga)
$$
is a smooth curve of projections (on the direct sum of all 
eigenspaces corresponding to eigenvalues in the interior of $\ga$) 
with finite dimensional ranges, the ranks
(i\.e\. dimension of the ranges) must be constant: it is easy to see 
that the  
(finite) rank cannot fall locally, and it cannot increase, since the 
distance in $L(H,H)$ of $P(t)$ to the subset of operators of 
rank $\le N=\operatorname{rank}(P(s))$ is continuous in $t$ and is either 
0 or 1.   
So for $t$ near $s$, there are equally many eigenvalues in the 
interior, and we may call them 
$\mu_i(t): 1\le i\le N$ (repeated 
with multiplicity).
Let us denote by 
$e_i(t): 1\le i\le N$ a corresponding system of eigenvectors of $A(t)$.
Then by the residue theorem we have 
$$
\sum_{i=1}^{N}\mu_i(t)^p e_i(t)\langle e_i(t),\quad \rangle
     = -\frac1{2\pi i}\int_\ga z^p(A(t)-z)\i\;dz 
$$
which is smooth in $t$ near $s$, as a curve of operators in $L(H,H)$ 
of rank $N$, since 0 is not an eigenvalue.
\roster
\item[2]
{\sl Claim.} Let $t\mapsto T(t)\in L(H,H)$ be a smooth curve of 
operators of rank $N$ in Hilbert space such that $T(0)T(0)(H)=T(0)(H)$.
Then $t\mapsto \operatorname{Trace}(T(t))$ is smooth near $0$
(note that this implies $T$ smooth into the space of nuclear 
operators, since all bounded linear functionals are of the form 
$A\mapsto \operatorname{Trace}(AB)$, by \cit!{9},~2.3 or 2.14.(4).
\endroster

Let $F:=T(0)(H)$.
Then $T(t)=(T_1(t),T_2(t)):H\to F\oplus F^\bot$ and 
the image of $T(t)$ is the space 
$$\align
T(t)(H) &= \{(T_1(t)(x),T_2(t)(x)):x\in H\}\\
&= \{(T_1(t)(x),T_2(t)(x)):x\in F\}\text{ for }t\text{ near }0\\
&= \{(y,S(t)(y)):y\in F\}, \text{ where }S(t):=T_2(t)\o (T_1(t)|F)\i.\\
\endalign$$
Note that $S(t):F\to F^\bot$ is smooth in $t$ by finite 
dimensional inversion for $T_1(t)|F:F\to F$.
Now
$$\align
\operatorname{Trace}(T(t)) &= \operatorname{Trace}\left(  
     \pmatrix 1 & 0 \\ -S(t) & 1 \endpmatrix 
     \pmatrix T_1(t)|F & T_1(t)|F^\bot \\ 
          T_2(t)|F & T_2(t)|F^\bot\endpmatrix
     \pmatrix 1 & 0 \\ S(t) & 1 \endpmatrix
     \right)\\
&= \operatorname{Trace}\left(  
     \pmatrix T_1(t)|F & T_1(t)|F^\bot \\ 0 
          & -S(t)T_1(t)|F^\bot+T_2(t)|F^\bot \endpmatrix
     \pmatrix 1 & 0 \\ S(t) & 1 \endpmatrix\right)\\
&= \operatorname{Trace}\left(  
     \pmatrix T_1(t)|F & T_1(t)|F^\bot \\ 0 & 0 \endpmatrix
     \pmatrix 1 & 0 \\ S(t) & 1 \endpmatrix\right),
     \text{ since rank}=N\\
&= \operatorname{Trace}
     \pmatrix T_1(t)|F+(T_1(t)|F^\bot)S(t)  & T_1(t)|F^\bot \\ 
          0 & 0 \endpmatrix\\
&= \operatorname{Trace}\Bigl(T_1(t)|F+(T_1(t)|F^\bot)S(t):F\to F\Bigr),
\endalign$$ 
which is visibly smooth since $F$ is finite dimensional. 

 From claim \therosteritem2 we now may conclude that 
$$
\sum_{i=1}^{N}\mu_i(t)^p =
-\frac1{2\pi i}\operatorname{Trace}\int_\ga z^p(A(t)-z)\i\;dz 
$$
is smooth for $t$ near $s$.

Thus the Newton polynomial mapping
$s^N(\mu_{1}(t),\dots,\mu_N(t))$ is smooth, so also the elementary 
symmetric polynomial 
$\si^N(\mu_{1}(t),\dots,\mu_N(t))$ is smooth, and thus  
$\{\mu_i(t): 1\le i \le N\}$ is the set of roots of a polynomial 
with smooth coefficients. 
By theorem \nmb!{4.3} there is an arrangement of these roots such 
that they become differentiable.
If no two of the generically 
different ones meet of infinite order, by theorem \nmb!{4.2} there is 
even a smooth arrangement. 

To see that in the general case they are even $C^1$ 
note that the images of the 
projections $P(t,\ga)$ of constant rank for $t$ near $s$ describe 
the fibers of a smooth vector bundle.
The restriction of $A(t)$ to 
this bundle, viewed in a smooth framing, becomes a smooth curve of 
symmetric matrices, for which by Rellich's result \nmb!{7.3} the 
eigenvalues can be chosen $C^1$. 
This finishes the proof of claim \therosteritem1.

\roster
\item[3]
{\sl Claim.} 
Let $t\mapsto\la_i(t)$ be a differentiable eigenvalue of $A(t)$, 
defined on some interval.
Then
$$
|\la_i(t_1)-\la_i(t_2)| \le (1+|\la_i(t_2)|)(e^{a|t_1-t_2|}-1) 
$$
holds for a continuous positive function $a=a(t_1,t_2)$ which is 
independent of the choice of the eigenvalue.
\endroster
For fixed $t$ near $s$ take all roots $\la_j$ which meet $\la_i$ at 
$t$, order them differentiably near $t$, and consider the projector 
$P(t,\ga)$ onto the joint eigenspaces for only those roots (where 
$\ga$ is a simple smooth curve containing only $\la_i(t)$ in its 
interior, of all the  
eigenvalues at $t$).
Then the image of $u\mapsto P(u,\ga)$, for $u$ 
near $t$, describes a smooth finite dimensional vector subbundle of 
$\Bbb R\x H$, since its rank is constant.
For each $u$ choose an 
othonormal system of eigenvectors $v_j(u)$ of $A(u)$ corresponding to 
these $\la_j(u)$.
They form a (not necessarily continuous) framing of 
this bundle. 
For any sequence $t_k\to t$ there is a subsequence such that each 
$v_j(t_k)\to w_j(t)$ where $w_j(t)$ is again an orthonormal system of 
eigenvectors of $A(t)$ for the eigenspace of $\la_i(t)$. 
Now consider 
$$
\frac{A(t)-\la_i(t)}{t_k-t}v_i(t_k) +
\frac{A(t_k)-A(t)}{t_k-t}v_i(t_k) -
\frac{\la_i(t_k)-\la_i(t)}{t_k-t}v_i(t_k) = 0,
$$
take the inner product of this with $w_i(t)$, note that then the 
first summand vanishes, and let $t_k\to t$ to obtain
$$
\la_i'(t) = \langle A'(t)w_i(t),w_i(t) \rangle
\text{ for an eigenvector }w_i(t)\text{ of }A(t)
\text{ with eigenvalue }\la_i(t).
$$
This implies, where $V_t=(V,\|\quad\|_t)$, 
$$\align
|\la_i'(t)|&\le \|A'(t)\|_{L(V_t,H)}
     \|w_i(t)\|_{V_t}\|w_i(t)\|_H \\
&= \|A'(t)\|_{L(V_t,H)}
     \sqrt{\|w_i(t)\|_H^2+\|A(t)w_i(t)\|_H^2} \\
&= \|A'(t)\|_{L(V_t,H)}\sqrt{1+\la_i(t)^2}\le a+a|\la_i(t)|,\\
\endalign$$
for a constant $a$ which is valid for a compact interval of $t$'s 
since $t\mapsto \|\quad\|_t^2$ is smooth on $V$.
By Gronwall's lemma (see e\.g\. \cit!{3}, (10.5.1.3)) this implies 
claim \therosteritem3.

By the following arguments we can conclude that all eigenvalues 
may be numbered as $\la_i(t)$ for $i$ in $\Bbb N$ or $\Bbb Z$  
in such a way that they are differentiable (by which we mean $C^1$, or 
$C^\infty$ under the stronger assumption) in $t\in \Bbb R$.
Note first that by claim \therosteritem3 no eigenvalue can go off to 
infinity in  
finite time since it may increase at most exponentially.
Let us first number all eigenvalues of $A(0)$ increasingly. 

We claim that for one eigenvalue (say $\la_0(0)$) there exists a 
differentiable extension to all of $\Bbb R$; namely the set of all 
$t\in \Bbb R$ with a differentiable extension of $\la_0$ on the 
segment from 0 to $t$ is open and closed.
Open follows from claim \therosteritem1.
If this intervall does not reach infinity,
from claim \therosteritem3 it follows that $(t,\la_0(t))$ has an 
accumulation  point $(s,x)$ at the the end $s$.
Clearly $x$ is an 
eigenvalue of $A(s)$, and by claim \therosteritem1 the eigenvalues 
passing through $(s,x)$ can be arranged differentiably, and thus 
$\la_0(t)$ converges to $x$ and can be extended differentiably beyond 
$s$. 

By the same argument we can extend iteratively all eigenvalues 
differentiably to all $t\in \Bbb R$: if it meets an already 
chosen one, the proof of \nmb!{4.3} shows that we may pass through it 
coherently.

Now we start to choose the eigenvectors smoothly, under the stronger 
assumption. 
Let us consider again eigenvalues $\{\la_i(t): 1\le i \le N\}$ 
contained in the interior of a smooth curve $\ga$ for $t$ in an 
open interval $I$.
Then $V_t:=P(t,\ga)(H)$ is the fiber of a smooth 
vector bundle of dimension $N$ over $I$.
We choose a smooth framing 
of this bundle, and use then the proof of theorem \nmb!{7.6} to 
choose smooth sub vector bundles whose fibers over $t$ are the 
eigenspaces of the eigenvalues with their generic multiplicity. 
By the same arguments as in \nmb!{7.6} we then get global vector sub 
bundles with fibers the eigenspaces of the eigenvalues with their 
generic multiplicity, and finally smooth eigenvectors for all 
eigenvalues.  
\qed\enddemo

\enddocument